# Computing Customers Sojourn Times in Jackson Networks Distribution Functions and Moments

Manuel Alberto M. Ferreira[1*]

## Abstract

Jackson queuing networks have a lot of practical applications, mainly in services and technologic devices. For the first case, an example are the healthcare networks and, for the second, the computation and telecommunications networks. Evidently the time that a customer - a person, a job, a message … – spends in this kind of systems, its sojourn time[1], is an important measure of its performance, among others. In this work, the practical statistical known results about the sojourn time of a customer, in a Jackson network, distribution are collected and presented. And an emphasis is set on the numerical methods applicable to compute the distribution function and the moments.

## 1 Introduction

The aim of this work is to present some problems and results that arise in the study of the customers sojourn times in Jackson networks of queues. These networks have many applications, namely in the modelling of healthcare, computation, and telecommunications networks. And a customer sojourn time, in this kind of system, is evidently an important element to be considered in its performance evaluation. Maybe the most important.

The network model to be considered in this paper is briefly described in section 2. Some essential results are there mentioned. The main objective of section 3 is the presentation of network flow equations**,** formula (3.2) that, in some situations allows the sojourn times moments exact computation. In section 4 it is given a numerical method for the sojourn times distribution function and any order moments computation, adequate to any Jackson network.

## 2 Results and Examples

Along this work, the sojourn times in a class of Markovian networks of queues, introduced initially by Jackson, see [1,2], will be studied. They are called Jackson networks and have only one class of customers. They are composed of $J$ nodes numbered $1,2,...,J$. It is usual to put $U = \{1,2, ... , J\}$.

In each node there is only one server, a queue discipline first-come-first-served" (*FCFS*) and an infinite waiting capacity.

They are open networks since any customer may enter or abandon it.
The exogenous arrivals-that is: from the outside of the network-process at node $j$ is a Poisson process (that is: the customers arrive one each time, and the interarrival times are independent and identically exponentially distributed) at rate $v_j, j \in U$, independent of the exogenous arrivals processes to the other nodes. It is stated that $v = \sum_{j=1}^{J} v_j$.

---

[1] The ***sojourn time*** *of a customer in a node is the sum of its **waiting time**, that is: the time the customer is in queue waiting to be served, plus its **service time**. Note, for curiosity, that in infinite servers' queues, which practical implementation is done guaranteeing that upon its arrival at the system a customer finds immediately an available server, the service time is the sojourn time.*

______________________________________________________________________________________________________________

[1]ISTAR-IUL - Information Sciences, Technologies and Architecture Research Center (ISTA), Instituto Universitário de Lisboa, Lisboa, Portugal.
*Corresponding author: E-mail: manuel.ferreira@iscte.pt

The service times at node $j$ are independent and identically distributed, having exponential distribution with parameter $\mu_j, j \in U$, and independent from the other nodes service times.

After the completion of a service at node $j$, a customer is immediately directed to node $l$ with probability $p_{jl}$, or abandons the network with probability $q_j = 1 - \sum_{l=1}^{J} p_{jl}, j \in U$. These probabilities are not influenced by the movements of the other customers in the network. The $p_{jl}$ matrix is called $P$. The matrix $P$ is called the commutation matrix and the $p_{jl}$ commutation probabilities.

The total arrivals rate, exogenous and endogenous-that is: from the other nodes of the network- at node $j$, $\theta_j$ satisfies the network traffic equations:

$$\theta_j = v_j + \sum_{l=1}^{J} \theta_l p_{lj}, j = 1,2, ..., J \tag{2.1}$$

The state of the network at instant $t$ is given by $N(t) = [N_1(t), ..., N_J(t)]$, where $N_j(t)$ is the number of customers at node $j$ in instant $t, j = 1,2, ..., J$. $N$ is the population process. If the traffic intensity $\rho_j = \frac{\theta_j}{\mu_j} < 1, j = 1,2, ..., J$ the process $N = \{N(t)\}$ possesses stationary, or equilibrium, distribution (when in equilibrium, the population process distribution does not change. In Jackson networks, this equilibrium distribution coincides with the distribution obtained by making $t$ converge to infinite, called the limit distribution), see for instance [3]:

$$\pi(n_1, n_2, ..., n_J) = \prod_{j=1}^{J} (1 - \rho_j) \rho_j^{n_j}, n_j \geq 0, j = 1,2, ..., J \tag{2.2}$$

The distribution (2.2) is of product form kind, see for instance [4,5], that is a very relevant concept in networks of queues population process equilibrium distributions.

Calling $S_j, W_j$ and $X_j$ the sojourn, waiting and service, respectively, times of a customer at node $j$

$$S_j = W_j + X_j \tag{2.3}$$

The Jackson networks sojourn times considered in this paper are those of typical customers that, arriving at the network, find the population process in an equilibrium state.

Call $S$ the sojourn time in the network, that is: the time that goes between the arrival at the network and the departure of one of those customers from it. If in its path it navigates the nodes $1,2, ..., l, S = S_1 + S_2 + +S_l$.

To study the sojourn time, the following notions are important:

- A network has feedback" if a customer may come back to the same node after the completion of its service, immediately or in a future moment,
- A network without feedback" is an acyclic" one,
- A network has "overtaking" if a customer can overtake" another one taking an alternative path between two nodes.

Now, three examples of typical Jackson networks usually considered in the study of sojourn times are presented. It may be said that more complex Jackson networks are integrated by networks fulfilling the properties of these examples, in a modular way.

## 2.1 Simple queues series

Indeed, this Jackson network is a series of $M|M|1$ queues. According to Kendal's notation here followed, the first $M$ means that customer arrivals at the node follow a Poisson process, the second $M$ that the length of the service time possesses an exponential distribution, and the 1 that there is only one server.

For this Jackson network

$$p_{jl} = \begin{cases} 1, if\ l = j+1, j = 1,2,\ldots,J-1 \\ 0, otherwise \end{cases}$$

$v_1 = v$, $v_j = 0, j = 2,\ldots,J$ and $\theta_j = v, j = 1,2,\ldots,J$. Fig. 1 is a graphical representation of a simple queue's series.

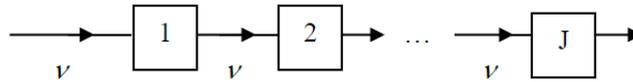

**Fig. 1. Simple queues series**

Some significant results for the simple queues series are:

- All customers' flows, in this network, at stationary state, are Poisson processes. It is a consequence of it, in stationary state, that the departure process from an $M|M|1$ queue is a Poisson process, see for instance [6],
- The sojourn times in the various nodes are independent random variables. In [6] it is presented a demonstration of this statement based on the reversibility concept,
- The sojourn time at node $j$ is an exponential random variable with parameter $\mu_j - v, j = 1,2,\ldots,J$. So, if $\mu_j = \mu, j = 1,2,\ldots,J$, the total sojourn time in this network is distributed as a $J$ order Erlang distribution with parameter $\mu - v$. This distribution function will be designated $E_{J,\mu-v}(t)$.
- The waiting times are dependent random variables. See also [6].

So, the sojourn time study in these networks has no difficulty. The same is not true for the waiting time.

## 2.2 $M|M|1$ queue with instantaneous Bernoulli feedback

It is a network with a single node. $J = 1$, $p_{11} = p$, $q_1 = 1 - p$ and $\theta = \frac{v}{1-p}$, where $\theta = \theta_1$ and $v = v_1$, see Fig. 2.

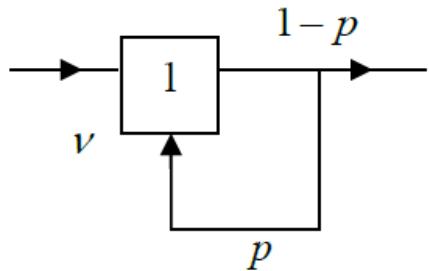

**Fig. 2. $M|M|1$ queue with Instantaneous Bernoulli Feedback**

It is an $M|M|1$ queue in which each customer, when leaving it after being served, returns to the queue with probability $p$.

Call $S_m$ the $m^{th}$ customer sojourn time in the network. So, if it is served $k$ times,

$$S_m = (t_{m1}^o - t_{m1}^a) + (t_{m2}^o - t_{m2}^i) + + (t_{mk-1}^o - t_{mk-1}^i) + (t_{mk}^d - t_{mk}^i) \quad (2.4),$$

Where

- $t_{ml}^0 - t_{ml}^i$ is the time that the customer spends passing by the service system in the *l*th time, given by the difference between the *l*th output (0) instant from the server and the one of the *l*th junction ( *i-input*) to the queue,
- $t_{m1}^0 - t_{m1}^a$ is the time that the customer spends passing by the service system for the first time, given by the difference between the first output (0) instant from the server and the one of the arrival (*a*) to the queue,
- $t_{mk}^d - t_{mk}^i$ is the time that the customer spends passing by the service system for the last time, given by the difference between the departure (*d*) instant from the network and the one of the $k^{th}$ junction (*i*) to the queue.

Note that $K$, the number of times that the customer passes by the server, is a random variable and $P(K = k) = (1-p)p^{k-1}, k = 1,2, \ldots$

The set $\{(t_{ml}^0 - t_{ml}^i): l = 2,3, \ldots\}$ is not a sequence of independent random variables, see [3]

So, it is not possible to make use of the usual statement to sum independent random variables. But it is possible to get an expression to $P(S_m \leq t)$ that requires the $k$ steps transition probabilities computation for the delayed Markovian renewal process:

$$\{(N(t^i - 0)(t_l^o - t_l^i)\}l = 0,1,2,\ldots\} \tag{2.5}$$

conditioning to the number of times that the customer returns to the queue. Calling that transition probabilities matrix $Q_i^k(t)$, see still [3],

$$P(S_m \leq t) = \sum_{k=1}^{\infty} \pi Q_i^k(t) p(1-p) V \tag{2.6}$$

where $\pi$ is the $N^i$ (embedded version of $N$ in the input instants) stationary distribution, $k$ is the number of times the customer passes by the server and $V$ is a vector which entries are all 1.

So, now, the situation is much more complicated than in the former case owing to the feedback.

## 2.3 The Jackson three node acyclic network

It is a network with three nodes each one behaving as an $M|M|1$ queue, Fig. 3, where $p_{12} = p$, $p_{13} = 1-p$, $p_{23} = 1$, $p_{jl} = 0$ in the other cases, $v_1 = v$, $v_j = 0, j = 2,3$, $\theta_1 = v$, $\theta_2 = pv$ and $\theta_3 = v$.

In equilibrium, all customers' flows are Poisson process in this network.

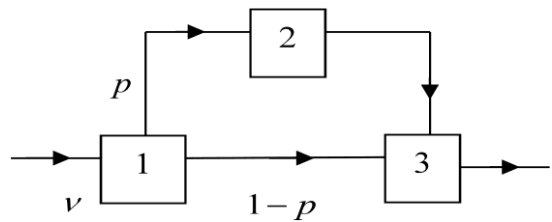

Fig. 3. Jackson three node acyclic network

Consequently,

- The sojourn time at node $j$ is a random variable exponentially distributed with parameter $\mu_j - \theta_j, j = 1,2,3$. $S_1$ and $S_2$ are independent random variables as well as $S_2$ and $S_3$.

This result is valid for any Jackson acyclic network:

- Suppose that a customer follows a path $r$ in a Jackson acyclic network with only one server at each node. If node $j$ belongs to path $r$, $S_j$ is such that fulfils the condition:

$$P(S_j \leq t | \text{the followed path is } r) = 1 - e^{-(\mu_j - \theta_j)t}, t \geq 0 \qquad (2.7)$$

and, if node $j$ is the next to the customer after node $l$, $S_j$ and $S_l$ are independent random variables.

But,

- $S_1$ and $S_3$ are dependent random variables: In [7] is stated that, indeed, $S_1$ and $S_3$ are positively correlated. In [8] it is stated that if

$$2(v + \mu_1 + \mu_2 + \mu_3)(\frac{1}{v + \mu_1 + \mu_2 + \mu_3}(\frac{1}{\mu_1 - v} + \frac{1}{\mu_3 - v}) + \frac{1}{(\mu_1 - v)^2}$$

$$+ \frac{1}{(\mu_3 - v)^2})^{\frac{1}{2}} < 1 \qquad (2.8)$$

and if

$$\frac{\mu_2 \left( \frac{1}{2(v+\mu_1+\mu_2+\mu_3)} - \frac{1}{2} \left( \frac{1}{(v+\mu_1+\mu_2+\mu_3)^2} + \frac{4}{v+\mu_1+\mu_2+\mu_3} \left( \frac{1}{\mu_1-v} + \frac{1}{\mu_3-v} \right) - 4 \left( \frac{1}{(\mu_1-v)^2} + \frac{1}{(\mu_3-v)^2} \right) \right)^{\frac{1}{2}} \right)}{1-v \left( \frac{1}{2(v+\mu_1+\mu_2+\mu_3)} - \frac{1}{2} \left( \frac{1}{(v+\mu_1+\mu_2+\mu_3)^2} + \frac{4}{v+\mu_1+\mu_2+\mu_3} \left( \frac{1}{\mu_1-v} + \frac{1}{\mu_3-v} \right) - 4 \left( \frac{1}{(\mu_1-v)^2} + \frac{1}{(\mu_3-v)^2} \right) \right)^{\frac{1}{2}} \right)} < p <$$

$$< \frac{\mu_2 \left( \frac{1}{2(v+\mu_1+\mu_2+\mu_3)} + \frac{1}{2} \left( \frac{1}{(v+\mu_1+\mu_2+\mu_3)^2} + \frac{4}{v+\mu_1+\mu_2+\mu_3} \left( \frac{1}{\mu_1-v} + \frac{1}{\mu_3-v} \right) - 4 \left( \frac{1}{(\mu_1-v)^2} + \frac{1}{(\mu_3-v)^2} \right) \right)^{\frac{1}{2}} \right)}{1-v \left( \frac{1}{2(v+\mu_1+\mu_2+\mu_3)} + \frac{1}{2} \left( \frac{1}{(v+\mu_1+\mu_2+\mu_3)^2} + \frac{4}{v+\mu_1+\mu_2+\mu_3} \left( \frac{1}{\mu_1-v} + \frac{1}{\mu_3-v} \right) - 4 \left( \frac{1}{(\mu_1-v)^2} + \frac{1}{(\mu_3-v)^2} \right) \right)^{\frac{1}{2}} \right)}$$

(2.9)

verify both simultaneously it is possible to guarantee that $S_1$ and $S_3$ are positively correlated in equilibrium.

Why this happen? One explanation may be the following:

- There are two alternative paths for a customer to go from node 1 to node 3. And a customer that follows by node 2 may be overtaken by another one that goes directly from node 1 to node 3. So, a customer, when arriving at node 3, may meet there another one that was behind it at node 1 or even that had not arrived when it was there.
- These overtaking customers can delay a certain customer, when it arrives at node 3, for a longer time than that if they were not present. The number of these customers depends, partly, on the number of the customers that arrive while the customer that is being followed is in node 1, partly owing to the supposition of a *FCFS* discipline.
- Consequently, the time that a customer waits at node 3 depends on how much time it has waited at node 1.

Now the complication is due to the overtaking.

For a more comprehensive treatment of this type of network, see [9]. For practical applications, see some examples in [10].

## 3 Network Flow Equations

The objective of this section is to present the so-called network flow equations" for the Jackson networks, that allow the deduction of formulae to the computation of sojourn times moments of any order, efficient in some situations.

Following the work of [11] call $\tau_j$ an arrival instant, endogenous or exogenous, at node $j$ and $\tau_j + T_j$ the departure instant from the network of the customer that arrived in $\tau_j, j = 1,2, j$, so

$T_j$ is the remaining sojourn time, in the network, for the arrival at node $j$ in the instant $\tau_j, j = 1,2,\ldots,J$.

Call $h_j$ the Laplace transform of the $T_j, j = 1,2, J$ distribution. As $N$ is a strong Markov process, and the network state process "seen by the arrivals" is in equilibrium, the $T_j, j = 1,2, J$ and its Laplace transforms are uniquely determined.

Dealing with the sojourn time as the lifetime of a Markov process $\vartheta$ it is possible to show that the Laplace transforms $h_j, j = 1,2, , J$ satisfy an equations system called the "network flow equations". That is, see [11],

Being $H_j$ the probability distribution with Laplace transform $h_j$, there is a distribution probability with Laplace transform $q_j$ such as

$h_j(s) + \frac{sg_j(s)}{\mu_j - \theta_j} = q_j(s) + \sum_{l=1}^{j} p_{jl} h_l(s) \geq 0$ and $j = 1,2, J$ (3.1).

In Jackson networks without overtaking" the transforms $h_j$ and $g_j$ are identical for each $j$. Given $h_j, j = 1,2, J$ the transforms $g_j, j = 1,2, J$ are uniquely determined by (3.1). The converse is also true since $I - P$, being $I$ the identity matrix, is invertible.

After (3.1), by successive derivations it is obtained:

### 3.1 Network flow equations

For $j = 1,2, J$ and $r = 1,2, \ldots$

$$E[T_j^r] = r! (\mu_j - \theta_j)^{-r} + \sum_{l=1}^{J} p_{jl} E[T_l^r]$$
$$+ \sum_{l=1}^{J} p_{jl} \sum_{n=1}^{r-1} \frac{r!}{n!(r-n)!} \mu_j^{-n} E[T_l^{r-n} \prod_{m=1}^{n}(N_j(T_l^-) + m)] \quad (3.2).$$

For $r = 1$, (3.2) assumes the matrix form

$$[E[T_j]] = (I - P)^{-1}[(\mu_j - \theta_j)^{-1}] \quad (3.3).$$

For $r = 2$, (3.2) assumes the form

$$E[T_j^2] = 2(\mu_j - \theta_j)^{-2}$$
$$+ \sum_{l=1}^{J} p_{jl} E[T_l^2] + 2\mu_j^{-1} \sum_{l=1}^{J} p_{jl} E[T_l(N_j(T_l^-) + 1)], j = 1,2, , J \quad (3.4).$$

**Note:**

1. Equality (3.4) defines a system of $J$ equations and $J^2 + J$ unknowns. In general, when $r \geq 2$, the product terms involving the variables $T_l$ and $N_j(\tau_l^-)$ prevent the exact computation of the sojourn times $r$ order moments; there are too many unknowns and too few equations. In these instances, other independent equations are needed to complement (3.4) to be possible to obtain exact solutions.

2. When any pair of nodes in the network is connected by, in the maximum, one oriented path and $p_{jj} = 0, j = 1,2,\ldots,J$, $T_l$ and $N_j(\tau_l^-)$ are independent for $j \neq l$. The computation of $E[T_j(N_j(\tau_j^-) + 1)]$ is irrelevant since $p_{jj} = 0, j = 1,2,\ldots,J$. In this case (3.2) becomes a compact recursive formula that allows the computation of any order moments of the sojourn times, $T_j, j = 1,2,\ldots J$. For instance, as, in these conditions,

$$E[N_j(\tau_l^-)] = \frac{\theta_j}{\mu_j - \theta_j}, j = 1,2,\ldots,J \quad (3.5),$$

(3.4) assumes the form:

$$E[T_j^2] = 2(\mu_j - \theta_j)^{-2} + \sum_{l=1}^{J} p_{jl} E[T_l^2] + 2(\mu_j - \theta_j)^{-1} \sum_{l=1}^{J} p_{jl} E[T_l] \; j = 1,2,J \quad (3.6).$$

Applying (3.6) to the simple queue's series:

$$E[T_j^2] = 2(\mu_j - v)^{-2} + E[T_{l+1}^2] + 2(\mu_j - v)^{-1} E[T_{l+1}], j = 1,2,\ldots,J1-1$$
$$E[T_J^2] = 2(\mu_J - v)^{-2}$$
$$E[T_J^2] = 2(\mu_J - v)^{-2} \quad (3.7)$$

that together with (3.3) results in

$$VAR[T_j] = \sum_{l=j}^{J} (\mu_l - v)^{-2} \quad (3.8).$$

3. *If* those conditions are not fulfilled, in [11] it is suggested to identify adequate Martingale families in $N$ as a process to determine independent equations to complement (3.2). Applying this proceeding to the $M|M|1$ queue with Instantaneous Bernoulli Feedback it was obtained:

$$VAR[T] = \frac{1}{((1-p)\mu - v)^2} \frac{(1-p^2)\mu + vp}{(1-p^2)\mu - vp} \quad (3.9)$$

and

$$COV[N(\tau^-), T] = \frac{v(1-p)\mu}{(1-p^2)\mu - vp} \quad (3.10). \blacksquare$$

To see more uses for Laplace transform, in queues context, see for instance [12].

## 4 Customers Sojourn Times in Jackson Networks Distribution Functions and Moments Numerical Computations

Now it is given a general method, which key is the procedure called "randomization procedure", to approximate "first passage times" distributions in direct time Markov processes, being the sojourn times in queue systems a particular case.

Call $\aleph = \{X(T): T \geq 0\}$ a regular Markov process, in continuous time with a countable states space $E$ and a bounded matrix infinitesimal generator $Q$. The elements of $Q$ are designated $Q(x,y), x, y \in E$ and $Q(x) = \sum_{y \in E - \{x\}} Q(x,y)$. $\psi(T)$ designates the $X(t)$ state probability vector:

$$\psi_t(x) = P\{X(T) = x\}, x \in E \tag{4.1}$$

$X$ models the evolution of a queue system during the sojourn of a given, "marked", customer in it.

The states $of E$ have two main components:

i) The queue system state,
ii) The "marked" customer position.

Be

- $A$, the states subset that describes the system till the departure of the marked" customer, and
- $B$, the states subset that describes the system aft er the departure of that customer.

Evidently

- $\{A, B\}$ is a partition of $E$,

- If $T$ is the time that the process $\aleph$ spends in $A$ till attaining $B$, for the first time, $T$ is precisely the sojourn time of the marked" customer in the network.

It is supposed that $\aleph$ will remain in $B$, with probability 1 aft er having attained it for the first time. In fact, as the evolution of the system aft er the departure of the "marked" customer is irrelevant, it may be supposed that $B$ is a closed set. That is, the process $\aleph$ cannot come back to $A$ after reaching $B$. The quantity of interest is the $T$ distribution function, $\tau(T)$. Note that

$$\tau(T) = P\{T \leq t\} = P\{X(T) \in B\} = 1 - P\{X(t) \in A\}, t \geq 0 \tag{4.2}$$

since the presented hypotheses guarantee that $\{T \leq t\} = \{X(T) \in B\}$.

After [2] it is concluded that

The problem of computing $\tau(t)$ is equivalent to the one of computing the transient distribution of $X(t)$ in $A$ computation.

So, it is necessary to compute the vector $\psi_t, t \geq 0$. Being $P_t, t \geq 0$, the $\aleph n$ transition matrix,

$$\psi_t = \psi_0 P_t, T \geq 0 \tag{4.3}$$

and

$$P_t = exp(Qt) = \sum_{i=0}^{\infty} \frac{T^i}{i!} Q^i, t \geq 0 \tag{4.4}$$

The randomization procedure consists in using in (4.4) an equivalent representation; see [13]:

$$P_\tau = exp(-\alpha t) exp(\alpha t(I + \frac{1}{\alpha}Q)) = exp(-\alpha t) \sum_{i=0}^{\infty} \frac{a^i t^i}{i!} R^i \tag{4.5}$$

where

$$R = I + \frac{1}{\alpha} Q \tag{4.6}$$

is called randomized matrix. $I$ is the identity matrix, and $\alpha$ is a positive upper bound for the whole $Q \chi \in E$.
Note that, see [14,15],

Although the equation (4.5) seems more complex than (4.4), it accomplishes indeed more favorable computational properties. The most important is that $R$ is a stochastic matrix while $Q$ is not. Consequently, the computation using (4.5) is stable, and using (4.4) is not.

The randomization procedure has an interesting probabilistic meaning, useful to determine bounds for $\tau(t)$. Indeed, being $R$ a stochastic matrix, it defines a discrete time Markov process:

$$\Im = \{Y_n : n = 0, 1, \dots\} \qquad (4.7)$$

if it is assumed $Y_0 = X(0)$. With this procedure, the relation between the processes $\aleph$ and $\Im$ is quite simple as it will be seen next.

Extend the discrete time process $\Im$ to a continuous time Markov process, such that:

i) The time intervals between jumps are exponential random variables $i.i.d.$ with mean $1/\alpha$,
ii) The jumps are commanded by $R$.

In [13] it is shown that the resulting process is precisely the original process $\aleph$; but when there is a sequence of jumps in $\Im$ from the state $\chi \in E$ to itself' this will be noticed in $\aleph$ as a long sojourn in state $x$.

So, the randomization procedure may be interpreted as a sowing in the process $\aleph$ with "fake" random jumps between the true jumps. The resulting process, designated by $\overline{\aleph}$, at which the "fake" jumps are visible, has the same probabilistic structure than $\aleph$ but with an advantage:

The sequence of the jump instants $in$ $\overline{\aleph}$, "fake" and "true", is now a Poisson Process. This is not, in general, the case of $\aleph$.

Note that $Y_n$ is the state of $\overline{\aleph}$ in the instant of the $n^{th}$ jump, "fake" or "true".

Suppose that $\overline{\aleph}$ reaches the set $B$ in its $n^{th}$ jump. Consequently the $\overline{\aleph}$ sojourn time, and so also the $\aleph$, in $A$ is the sum of n exponential independent random variables with mean $1/\alpha$. That is, the sojourn time has a $n$ order Erlang distribution with parameter $\alpha$. Its distribution function will be designated $E_{n,\alpha}(t)$.

Be $h(n)$ the probability that $\overline{\aleph}$ reaches $B$ in its $n^{th}$ jump. Call $(\phi_n$ the state probability vector of $Y_n$:

$\phi_n = \psi_0 R^n$ (4.8).

The quantities $h(n)$ are given by the equivalent formulae:

$$h(n) = \begin{cases} \sum_{x \in B} \phi_0(x), n = 0 \\ \sum_{x \in A} \sum_{y \in B} (\phi_{n-1}(x) R(x,y), n > 0 \quad (4.9) \end{cases}$$

or

$$h(n) = \begin{cases} 1 - \sum_{x \in A} \phi_0(x), n = 0 \\ \sum_{x \in A} (\phi_{n-1}(x) - \sum_{x \in A} (\phi_n(x), n > 0 \quad (4.10). \end{cases}$$

Given the probabilities $h(n)$ and, noting that $\sum_{n=0}^{\infty} h(n) = 1$, it is obtained the following expression:

$$\tau(T) = \sum_{n=0}^{\infty} h(n) E_{n,\alpha}(t), \, t \geq 0 \qquad (4.11),$$

$$E[T^m] = \frac{1}{\alpha^m} \sum_{n=0}^{\infty} n(n+1)\ldots(n+m-1)h(n), \, m = 1, 2, \ldots \qquad (4.12).$$

The formula (4.12) for $m = 1$ is

$$E[T] = \frac{1}{\alpha} E[H] \qquad (4.13)$$

being $H$ the number of $\aleph$ jumps till reaching $B$. Expression (4.13) is Little's Formula in this queue's context.

Equation (4.12) allows obtaining simple bounds for $\tau(t)$ that may, in principle, to become arbitrarily close. Equation (4.12) allows obtaining a lower bound for $E[T^k]$, in principle, so close of $E[T^k]$ as wished. So, given any integer $k \geq 0$

$$L_k(t) \leq \tau(t) \leq U_k(t) \qquad (4.14)$$

where

$$L_k(t) = \sum_{n=0}^{k} h(n) E_{n,\alpha}(t), \, t \geq 0 \qquad (4.15),$$

$$U_k(t) = 1 - \sum_{n=0}^{k} h(n) \bar{E}_{n=0}(t), \, t \geq 0 \qquad (4.16)$$

and

$$E[T^m]_{L,k} \leq E[T^m], m = 1,2, . \qquad (4.17)$$

where

$$E[T^m]_{L,k} = \frac{1}{\alpha^m} \sum_{n=0}^{k} n(n+1)\ldots(n+m-1)h(n), m = 1,2, \qquad (4.18).$$

It is easy to prove the proposition:

## 4.1 Proposition

If, for any $\varepsilon > 0$, $k$ is chosen in accordance with the rule:

$$k = \min\{n \geq 0: \sum_{i=0}^{n} h(i) \geq 1 - \varepsilon\} = k(\varepsilon), \qquad (4.19)$$

or equivalently

$$J = \min\{n \geq 0: \sum_{x \in A} (\phi_n(x) \leq \varepsilon\} = J(\varepsilon) \qquad (4.20)$$

$|\tau(T) - L_{J(\varepsilon)}| \leq \varepsilon$ and $|\tau(T) - U_{J(\varepsilon)}| \leq \varepsilon$, uniformely in $\tau \geq 0$.

**Note:**

-The main problem in the application of the method presented, stays in the difficulty of the $h(n)$ computation. Indeed, for it, it is necessary to compute the vectors $(\phi_n$ but only in the subset $A$ ofthe state's space. When states space $E$ is finite, as it happens in the case of closed networks, both $h(n)$ and $(\phi_n$ can, at first glance, be computed exactly, apart the mistakes brought by the approximations.

In practice the states space is often infinite or, although finite, prohibitively great. In these situations, it is mandatory to truncate $E$. So, it must be considered a new level of approximation since the $h(n)$, ($\phi_n$, etc. must also be approximated now.

Indeed, what are viable to obtain is $h(n)$ lower bounds because the $E$ truncation is translated in probability loss [15]. So, with these $h(n)$ approximate values, (4.12) and (4.15) go on being valid but

- The uniform convergence property seen above is lost,
- The rules analogous to (4.17) and (4.18) are not equivalent. The one generated by (4.17) may be even unviable and in practice it is used only the one generated by (4.18), see [15]. ∎

Using this method, in [16] is shown that, in a Jackson three node acyclic network, the total sojourn time distribution function for a customer that follows the path integrated by the nodes 1, 2, and 3 is not the same obtained considering that $S_1, S_2$ and $S_3$ are independent although this one, designated by $F(t)$ is a good" approximation of that one. They show that in some cases it was not true the following:

$$F^L(t) \leq F(t) \leq F^U(t) \geq 0 \qquad (4.21)$$

being $F^L(t)$ and $F^U(t)$ the lower bound and the upper bound, respectively, of that customer sojourn time distribution function, obtained through the described method.

This conclusion is important because, despite the dependence between $S_1$ and $S_3$, $F(t)$ could be the $S$ distribution function. In [17], it is presented an example of dependent random variables which sum has the same distribution as if the random variables were independent.

Finally note that the formula (4.12), see [18], seems to be more efficient than (3.2), although only allows to obtain moments lower bounds, because its field of application is much greater.

## 5 Conclusions

The sojourn time has an evident practice interest. And is and has been intensively studied. Evidently the problem of the computation of the sojourn times in networks of queues is one of the most difficult in these networks study. Indeed, analytic solutions are the exception and not the rule. And, when existing, are quite rough.

Most of the known works only present results on sojourn time distributions for only one customer in paths without overtaking with *FCFS* disciplines in the nodes. It seems that still there are not results for simultaneous distributions of various customers sojourn times.

It follows, from the examples seen in section 2, that the sojourn times, at Jackson networks computations, difficulties occur when there are feedback and overtaking. In the first case the input server process is not a Poisson process, becoming everything more complex. In the second case dependencies exist among a customer sojourn times in the various nodes, simultaneously complicated and subtle, that make the total sojourn time computation difficult even if the sojourn times in each node are easy to compute.

From all this it results the interest of the methods presented in sections 3 and 4 to compute exactly and approximately the quantities related with the Jackson networks sojourn times. But, in these cases, the technical difficulties resulting from both the analytical and computational implementation of the calculations, which are very demanding from this point of view, must be considered.

## Acknowledgements

This work is financed by national funds through FCT −Fundação para a Ciência e Tecnologia, I.P., under the project $UIDB/04466/2021$. Furthermore, the author thanks the Instituto Universitário de Lisboa and ISTAR-IUL, for their support.

# Competing Interests

Author has declared that no competing interests exist.

**Biography of author(s)**

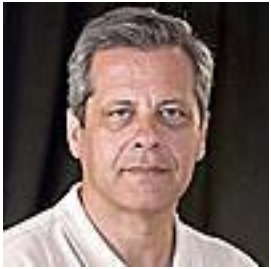

**Manuel Alberto M. Ferreira**
ISTAR-IUL - Information Sciences, Technologies and Architecture Research Center (ISTA), Instituto Universitário de Lisboa, Lisboa, Portugal

He is an Electrical Engineer and has a master's degree in Applied Mathematics by Lisbon University. He achieved his PhD in Management-Quantitative Methods and Habilitation in Quantitative Methods by ISCTE-Lisbon University Institute. He has been the former Chairman of Board of Directors and Vice-President of ISCTE-Lisbon University Institute, Full Professor (retired) of ISCTE-Lisbon University Institute. He is also the Emeritus Full Professor of ISCTE-Lisbon University Institute, former Director of Department of Mathematics in ISTA-School of Technology and Architecture, member of ISTAR-IUL (Integrated Researcher) and BRU-IUL (Associate Researcher) research centres, former member of UNIDE (Integrated Researcher) research centre. His research interests include Mathematics; Statistics; Stochastic Processes-Queues and Applied Probabilities; Game Theory; Chaos Theory; Bayesian Statistics: Application to Forensic Identification; Applications to Economics, Management, Business, Marketing, Finance, and Social Problems in general. He has published 468 papers, in Scientific Journals and Conference Proceedings, 43 book chapters, and 40 books and has presented 180 communications in International Scientific Conferences. He has 53 papers indexed in Scopus, 23 in Web of Science, 10 in zbMATH, and 66 in Crossref. He is the editor of several Scientific Journals and Conference Proceedings. He also reviewed papers for many Scientific Journals and Conferences.